\documentclass[12pt]{article}
\usepackage{mathrsfs}
\usepackage{amscd}
\usepackage{amsmath,amsfonts,amssymb,amscd}
\usepackage{indentfirst,graphics,epsfig,psfrag}
\input{epsf}
\usepackage{ifpdf}
\usepackage{enumerate}
\usepackage{appendix}
\usepackage{enumerate}

\usepackage{lineno}

\makeatletter \@addtoreset{figure}{section} \makeatother
\makeatletter
\long\def\@makecaption#1#2{%
   \vskip 10\p@
   \setbox\@tempboxa\hbox{{#1}\ \ #2}%
   \ifdim \wd\@tempboxa >\hsize
       {#1}\ \ #2\par
   \else
       \hbox to\hsize{\hfil\box\@tempboxa\hfil}%
   \fi}
\makeatother

\newtheorem{thm}{Theorem}[section]
\newtheorem{cor}[thm]{Corollary}
\newtheorem{lem}[thm]{Lemma}

\newtheorem{exe}[thm]{Example}
\newtheorem{obs}[thm]{Observation}

\newcommand{\qed}{{\hfill\rule{3pt}{7pt}}}

\def\qed{\hfill \rule{4pt}{7pt}}

\begin{document}
\title{\bf Upper bounds for the rainbow connection numbers of line graphs\footnote{Supported by
NSFC, PCSIRT and the ``973" program.
}}
\author{
\small  Xueliang Li,  Yuefang Sun\\
\small Center for Combinatorics and LPMC-TJKLC\\
\small Nankai University, Tianjin 300071, P.R. China\\
\small E-mail: lxl@nankai.edu.cn; syf@cfc.nankai.edu.cn
 }
\date{}
\maketitle
\begin{abstract}
A path in an edge-colored graph $G$, where adjacent edges may be
colored the same, is called a rainbow path if no two edges of it are
colored the same. A nontrivial connected graph $G$ is rainbow
connected if for any two vertices of $G$ there is a rainbow path
connecting them. The rainbow connection number of $G$, denoted by 
$rc(G)$, is defined as the smallest number of colors by using which
there is a coloring such that $G$ is rainbow connected. In this
paper, we mainly study the rainbow connection number of the line
graph of a graph which contains triangles and get two sharp upper
bounds for $rc(L(G))$, in terms of the number of edge-disjoint
triangles of $G$ where $L(G)$ is the line graph of $G$. We also give
results on the iterated line graphs.\\[2mm]
\end{abstract}

\section{Introduction}

All graphs in this paper are simple, finite and undirected. Let $G$
be a nontrivial connected graph with an edge coloring $c:
E(G)\rightarrow \{1,2,\cdots,k\}$, $k\in \mathbb{N}$, where adjacent
edges may be colored the same. A path of $G$ is called $rainbow$ if
no two edges of it are colored the same. An edge-colored graph $G$
is $rainbow ~connected$ if for any two vertices there is a rainbow
path connecting them. Clearly, if a graph is rainbow connected, it
must be connected. Conversely, any connected graph has a trivial
edge coloring that makes it rainbow connected, i.e., the coloring
such that each edge has a distinct color. Thus, we define the
$rainbow~connection~number$ of a connected graph $G$, denoted by
$rc(G)$, as the smallest number of colors for which there is an edge
coloring of $G$ such that $G$ is rainbow connected. An easy
observation is that if $G$ has $n$ vertices then $rc(G)\leq n-1$,
since one may color the edges of a spanning tree with distinct
colors, and color the remaining edges with one of the colors already
used. Generally, if $G_1$ is a connected spanning subgraph of $G$,
then $rc(G)\leq rc(G_1)$. We notice the trivial fact that $rc(G)=1$
if and only if $G$ is complete, and the fact that $rc(G)=n-1$ if and
only if $G$ is a tree, as well as the easy observation that a cycle
with $k>3$ vertices has rainbow connection number $\lceil
\frac{k}{2}\rceil$ (\cite{Chartrand 1}). Since a Hamiltonian graph
$G$ has a Hamiltonian cycle which contains all $n$ vertices, then
$G$ has rainbow connection number at most $\lceil \frac{n}{2}
\rceil$. Also notice that, clearly, $rc(G)\geq diam(G)$ where
$diam(G)$ denotes the diameter of $G$.

Chartrand et al. in \cite{Chartrand 1} determined that the rainbow
connection numbers of some graphs including trees, cycles, wheels,
complete bipartite graphs and complete multipartite graphs. Caro et
al. \cite{Y. Caro} gave some results on general graphs in terms of
some graph parameters, such as the order or the minimum degree of a
graph. They observed that $rc(G)$ can be bounded by a function of
$\delta(G)$, the minimum degree of $G$. They proved that if
$\delta(G)\geq 3$ then $rc(G)\leq \alpha n$ where $\alpha < 1$ is a
constant and $n=|V(G)|$. They conjectured that $\alpha=3/4$ suffices
and proved that $\alpha < 5/6$. Specifically, it was proved in
\cite{Y. Caro} that if $\delta=\delta(G)$ then $rc(G)\leq
\min\{\frac{\ln{\delta}}{\delta}n(1+o_{\delta}(1)),n\frac{4\ln{\delta}+3}{\delta}\}$.
Some special graph classes, such as line graphs, have many special
properties, and by these properties, we can get some interesting
results on their rainbow connection numbers in terms of some graph
parameters. For example, in \cite{Y. Caro} the authors got a very
good upper bound for the rainbow connection number of a 2-connected
graph according to their ear-decomposition. And in \cite{Li-Sun1},
we studied the rainbow connection numbers of line graphs of
triangle-free graphs in the light of particular properties of line
graphs of triangle-free graphs shown in \cite{S.T. Hedetniemi}, and
particularly, of 2-connected triangle-free graphs according to their
ear decompositions. However, we did not get bounds of the rainbow
connection numbers for line graphs that do contain triangles. In
this paper, we aim to investigate the remaining case, i.e., line
graphs that do contain triangles, and give two sharp upper bounds in
terms of the number of edge-disjoint triangles.

We use~$V(G)$, $E(G)$ for the sets of vertices and edges of $G$,
respectively. For any subset $X$ of $V(G)$, let $G[X]$ denote the
subgraph induced by $X$, and $E[X]$ the edge set of $G[X]$;
similarly, for any subset $E_1$ of $E(G)$, let $G[E_1]$ denote the
subgraph induced by $E_1$. Let $\mathcal {G}$ be a set of graphs,
then $V(\mathcal {G})=\bigcup_{G\in \mathcal {G}}{V(G)}$,
$E(\mathcal {G})=\bigcup_{G\in \mathcal {G}}{E(G)}$. We define a
$clique$ in a graph $G$ to be a complete subgraph of $G$, and a
$maximal~clique$ is a clique that is not contained in any larger
clique. The $clique~graph$ $K(G)$ of $G$ is the intersection graph
of the maximal cliques of $G$--that is, the vertices of $K(G)$
correspond to the maximal cliques of $G$, and two of these vertices
are joined by an edge if and only if the corresponding maximal
cliques intersect. Let $[n]=\{1,\cdots,n\}$ denote the set of the
first $n$ natural numbers. For a set $S$, $|S|$ denotes the
cardinality of $S$. We follow the notations and terminology of
\cite{Bondy} for those not defined here.

\section{Some basic observations}

We first list two observations which were given in \cite{Li-Sun1}
and will be used in the sequel.

\begin{obs}\label{bslem2.2}(\cite{Li-Sun1})
If $G$ is a connected graph and $\{E_i\}_{i\in [t]}$ is a partition
of the edge set of $G$ into connected subgraphs $G_i=G[E_i]$ and
$rc(G_i)=c_i$, then $$rc(G)\leq \sum_{i=1}^{t}{c_i}.$$
\end{obs} \qed

Let $G$ be a connected graph, and $X$ a proper subset of $V(G)$. To
$shrink$ $X$ is to delete all the edges between vertices of $X$ and
then identify the vertices of $X$ into a single vertex, namely $w$.
We denote the resulting graph by $G/X$.

\begin{obs}\label{bsthm3.6}(\cite{Li-Sun1}) Let $G'$ and $G$ be two
connected graphs, where $G'$ is obtained from $G$ by shrinking a
proper subset $X$ of $V(G)$, that is, $G'=G/X$, such that any two
vertices of $X$ have no common adjacent vertex in $V\setminus X$.
Then $$rc(G')\leq rc(G).$$ \qed
\end{obs}

Now we introduce two graph operations and two corresponding results
which will be used later.

\begin{figure}[!hbpt]
\begin{center}
\includegraphics[scale=1.0000000]{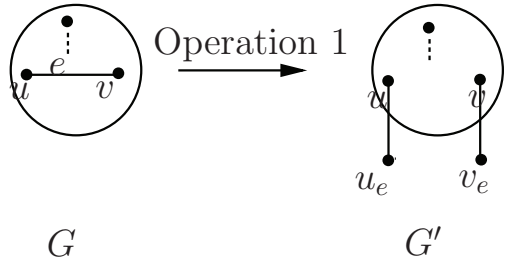}
\end{center}
\caption{$G'$ is obtained from $G$ by doing Operation 1 to edge
$e$.}\label{figure2}
\end{figure}

\noindent \textbf{Operation 1:} As shown in Figure \ref{figure2},
for any edge $e=uv \in G$ with $min\{deg_{G}(u),$ $deg_{G}(v)\}\geq
2$, we first subdivide $e$, then replace the new vertex with two new
vertices $u_e, v_e$ with $deg_{G'}(u_e)=deg_{G'}(v_e)=1$ where $G'$
is the new graph.

Since $deg_{G'}(u_e)=deg_{G'}(v_e)=1$ in $G'$, and by the definition
of a line graph, it is easy to show that $L(G)$ can be obtained from
$L(G')$ by shrinking a vertex set of two nonadjacent vertices (these
two vertices correspond to edges $uu_e$, $vv_e$, and belong to
cliques $\langle S(u)\rangle$, $\langle S(v)\rangle$, respectively,
in $L(G')$). Recall that the $line~graph$ of a graph $G$ is the
graph $L(G)$ whose vertex set $V(L(G))=E(G)$ and two vertices $e_1$,
$e_2$ of $L(G)$ are adjacent if and only if they are adjacent in
$G$. So by Observation \ref{bsthm3.6}, we have

\begin{obs}\label{bsthm2} If graph $G'$ is obtained from a connected
graph $G$ by doing Operation 1 at some edge $e \in G$, then
$$rc(L(G))\leq rc(L(G')).$$ \qed
\end{obs}

\begin{figure}[!hbpt]
\begin{center}
\includegraphics[scale=1.0000000]{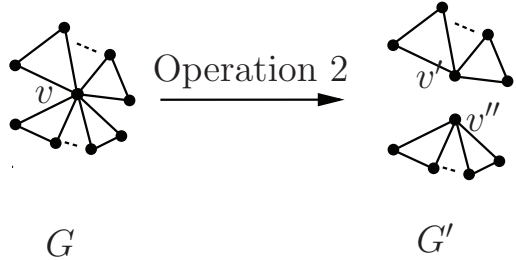}
\end{center}
\caption{Figure of Operation 2.}\label{figure6}
\end{figure}

\noindent \textbf{Operation 2}. As shown in Figure \ref{figure6},
$v$ is a common vertex of a set of edge-disjoint triangles in $G$.
We replace $v$ by two nonadjacent vertices $v'$ and $v''$ such that
$v'$ is the common vertex of some triangles, and $v''$ is the common
vertex of the rest triangles.

Since during this procedure, the number of edges does not change,
the order of the line graph $L(G)$ is equal to that of $L(G')$.
Furthermore, by the definition of a line graph, $L(G')$ is a
spanning subgraph of $L(G)$. So we have

\begin{obs}\label{obs2} If a connected graph $G'$ is obtained from a
connected graph $G$ by doing Operation 2 at some vertex $v\in G$,
then $$rc(L(G))\leq rc(L(G')).$$ \qed
\end{obs}

\section{Main results}

\subsection{A sharp upper bound}

Recall that the star, $S(v)$, at a vertex $v$ of graph $G$, is the
set of all edges incident to $v$. A $clique~decomposition$ of $G$ is
a collection $\mathscr{C}$ of cliques such that each edge of $G$
occurs in exactly one clique in $\mathscr{C}$.

We now introduce a new terminology. For a connected graph $G$, we
call $G$ a $clique$-$tree$-$structure$, if it satisfies the
following condition:

\textbf{$T_1$.} Each block is a maximal clique.

\begin{figure}[!hbpt]
\begin{center}
\includegraphics[scale=1.0000000]{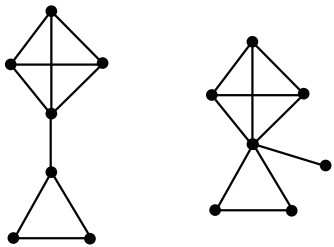}
\end{center}
\caption{A clique-forest-structure with size 6 and 2
components.}\label{figure1}
\end{figure}

We call a graph $H$ a $clique$-$forest$-$structure$, if $H$ is a
disjoint union of some clique-tree-structures, that is, each
component of a clique-forest-structure is a clique-tree-structure.
By condition $T_1$, we know that any two maximal cliques of $G$ have
at most one common vertex. Furthermore, $G$ is formed by its maximal
cliques. The $size$ of the clique-tree(forest)-structure is the
number of its maximal cliques. An example of clique-forest-structure
is shown in Figure \ref{figure1}.

If each block of a clique-tree-structure is a triangle, we call it a
$triangle$-$tree$-$structure$. Let $\ell$ be the size of a
triangle-tree-structure. Then, by definition, it is easy to show
that there are $2\ell+1$ vertices in it. Similarly, we can give the
definition of $triangle$-$forest$-$structure$. A
clique-tree-structure $G$ is called a $clique$-$path$-$structure$ if
the clique graph $K(G)$ is a path.

For a connected graph $G$, we call $G$ a
$clique$-$cycle$-$structure$, if it satisfies the following three
conditions:

\textbf{$C_1$.} $G$ has at least three maximal cliques;

\textbf{$C_2$.} Each edge belongs to exactly one maximal clique;

\begin{figure}[!hbpt]
\begin{center}
\includegraphics[scale=1.0000000]{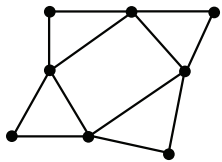}
\end{center}
\caption{An example of triangle-cycle-structure.}\label{figure11}
\end{figure}

\textbf{$C_3$.} The clique graph is a cycle. (In particular, if each
maximal clique is a triangle, then it is a
$triangle$-$cycle$-$structure$. An example of
triangle-cycle-structure is shown in Figure \ref{figure11}.)

An $inner~vertex$ of a graph is a vertex with degree at least two.
For a graph $G$, we use $\overline{V_2}$ to denote the set of all
inner vertices of $G$. Let $n_1=|\{v: deg_G(v)=1 \}|$,
$n_2=|\overline{V_2}|$. $\langle S(v)\rangle$ is the subgraph of
$L(G)$ induced by $S(v)$, clearly it is a clique of $L(G)$. Let
$\mathscr{K}_0=\{\langle S(v)\rangle: v\in V(G) \}$,
$\mathscr{K}=\{\langle S(v)\rangle: v\in \overline{V_2} \}$. It is
easy to show that $\mathscr{K}_0$ is a clique decomposition of
$L(G)$ (\cite{Bo-Jr Li}) and each vertex of the line graph belongs
to at most two elements of $\mathscr{K}_0$. We know that each
element $\langle S(v)\rangle$ of ${\mathscr{K}_0}\setminus
{\mathscr{K}}$, a single vertex of $L(G)$, is contained in the
clique induced by $u$ that is adjacent to $v$ in $G$. So
$\mathscr{K}$ is a clique decomposition of $L(G)$.

Let $X$ and $Y$ be sets of vertices of a graph $G$. We denote by
$E[X,Y]$ the set of the edges of $G$ with one end in $X$ and the
other end in $Y$. If $Y=X$, we simply write $E(X)$ for $E[X,X]$.
When $Y=V\backslash X$, the set $E[X,Y]$ is called the $edge~cut$ of
$G$ associated with $X$, and is denoted by $\partial(X)$.

\begin{thm}\label{bsthm3}
For any set $\mathcal {T}$ of $t$ edge-disjoint triangles of a
connected graph $G$, if the subgraph induced by the edge set
$E(\mathcal {T})$ is a triangle-forest-structure, then
$$rc(L(G))\leq n_2-t.$$ Moreover, the bound is sharp.
\end{thm}
\begin{pf}
Let $\mathcal {T}=\bigcup_{i=1}^{c}{\mathcal
{T}_i}=\bigcup_{i=1}^{c}\{T_{i,j_i}~is~ a~ triangle~ of~ $G$: 1\leq
j_i \leq t_i \}$($\sum_{i=1}^c{t_i}=t$) be a set of $t$
edge-disjoint triangles of $G$ such that the subgraph of $G$,
$G[E(\mathcal {T}_i)]$, induced by each $E(\mathcal {T}_i)$ is a
component of the subgraph $G[E(\mathcal {T})]$, that is, a
triangle-tree-structure of size $t_i$.

In $G$, for each $1\leq i\leq c$, let $G_i=G[E(\mathcal {T}_i)]$,
$V_i=V(G_i)$, $E_i=E(\mathcal {T}_i)$; $E_i^0=E(V_i)\cup
\partial(V_i)\supseteq E_i$, and $G_i^0=G[E_i^0]$. We obtain a new
graph $G'$ from $G$ by doing Operation 1 at each edge $e\in
E(V_i)\backslash E_i$ for $1\leq i\leq c$, and we denote by $G_i'$
the new subgraph (of $G'$) corresponding to $G_i^0$. Applying
Observation \ref{bsthm2} repeatedly, we have $rc(L(G))\leq
rc(L(G'))$.

Next we will show $rc(L(G'))\leq n_2-t$. By previous discussion, we
know that
$$\mathscr{K}=\{\langle S(v)\rangle: v\in \overline{V_2}
\}=\bigcup_{i=1}^c\{\langle S(v)\rangle: v\in V_i\}\bigcup \{\langle
S(v)\rangle: v\in \overline{V_2}\backslash \bigcup_{i=1}^c{V_i} \}$$
is a clique partition of $L(G)$. So
$$\{E(\langle S(v)\rangle): v\in V_i\}_{i=1}^c \bigcup \{E(\langle
S(v)\rangle): v\in \overline{V_2}\backslash \bigcup_{i=1}^c{V_i}
\},$$ that is, $$\{E(L(G_i^0))\}_{i=1}^c \bigcup \{E(\langle
S(v)\rangle): v\in \overline{V_2}\backslash \bigcup_{i=1}^c{V_i}
\}$$ is an edge partition of $L(G)$. So
$$\{E(L(G_i'))\}_{i=1}^c \bigcup \{E(\langle S(v)\rangle): v\in
\overline{V_2}\backslash \bigcup_{i=1}^c{V_i} \}$$ is an edge
partition of $L(G')$. By Observation \ref{bslem2.2}, we have
$$rc(L(G'))\leq \sum_{i=1}^c{rc(L(G_i'))}+\sum_{v\in
\overline{V_2}\backslash \bigcup_{i=1}^c{V_i}}{rc(\langle
S(v)\rangle)}.$$ We know $|V_i|=2t_i+1$, since the
triangle-tree-structure $G_i$ has size $t_i$. So $$rc(L(G'))\leq
\sum_{i=1}^c{rc(L(G_i'))}+(n_2-2t-c).$$ In order to get
$rc(L(G'))\leq n_2-t$, we need to show $rc(L(G_i'))\leq t_i+1$.

\textbf{Claim}. $rc(L(G_i'))\leq t_i+1$.

\textbf{Proof of the Claim}. Since the graph $G_i'$ is obtained from
$G_i$ by doing Operation 1 at each edge $e\in E(V_i)\backslash E_i$,
$G_i'$ contains exactly $t_i$ triangles: $\{T_{i,j_i}\in \mathcal
{T}_i: 1\leq j_i \leq t_i \}$. We will show that there is a
$(t_i+1)$-rainbow coloring of $L(G_i')$ by induction on $t_i$. For
$t_i=1$, $G_i'$ contains exactly one triangle, and we give its line
graph a 2-rainbow coloring as shown in Figure \ref{figure3}. We give
color 1 to the edges of $\langle S(u)\rangle$ incident with vertex
$e_1$, edges of $\langle S(v)\rangle$ incident with vertex $e_2$,
and edges of $\langle S(w)\rangle$ incident with vertex $e_3$; We
then give color 2 to the edges of $\langle S(u)\rangle$ incident
with vertex $e_3$, edges of $\langle S(v)\rangle$ incident with
vertex $e_1$, and edges of $\langle S(w)\rangle$ incident with
vertex $e_2$; Finally, give color 2 to the rest of the edges. It is
easy to show that this is a rainbow coloring. So, the above
conclusion holds for the case $t_i=1$.
\begin{figure}[!hbpt]
\begin{center}
\includegraphics[scale=1.0000000]{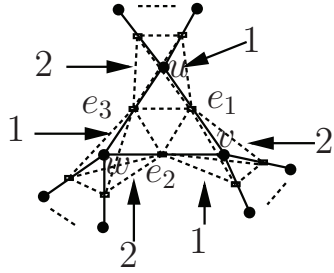}
\end{center}
\caption{2-rainbow coloring of line graph of graph with exactly one
triangle.}\label{figure3}
\end{figure}
We assume that the conclusion holds for the case $t_i=h-1$, and now
show that it holds for the case $t_i=h$. By the definition of
triangle-tree-structure, there must exist one triangle:
$T=\{u,v,w\}$, which has exactly one common vertex, namely $u$, with
other triangles and $v,w$ do not belong to any other triangle. We
now obtain a new graph $\overline{G_i'}$ from $G_i'$ by doing
Operation 1 at edges $e_1$ and $e_2$. It is clear that
$L(\overline{G_i'})$ can be obtained from $L(G_i')$ by subdividing
vertex $e_1$ into two vertices $\{e_1',e_1''\}$ and $e_2$ into two
vertices $\{e_2',e_2''\}$. Since $H_1$ has $h-1$ (edge-disjoint)
triangles, by induction hypothesis, $rc(L(H_1))\leq h$. So the
subgraph $L(G_i')\backslash \{\langle S(v)\rangle,\langle
S(w)\rangle\}\cong L(H_1)$ has a rainbow $h$-coloring, and we now
color the edges of $\langle S(v)\rangle$ and $\langle S(w)\rangle$
in the graph $L(G_i')$ as follows: give a new color to the edges of
$\langle S(w)\rangle$ incident with vertex $e_1$, and edges of
$\langle S(v)\rangle$ incident with vertex $e_2$. Let $e_3=vw$, we
then give any one color, say $c_1$, of the former $h$ colors to the
edges of $\langle S(w)\rangle$ incident with vertex $e_3$, and give
a distinct color, say $c_2$, of the former $h$ colors to the edges
of $\langle S(v)\rangle$ incident with vertex $e_3$. It is easy to
show that, with above coloring, $L(G_i')$ is rainbow connected, and
so $L(G_i')\leq t_i+1$ holds for $t_i=h$.
\begin{figure}[!hbpt]
\begin{center}
\includegraphics[scale=0.9500000]{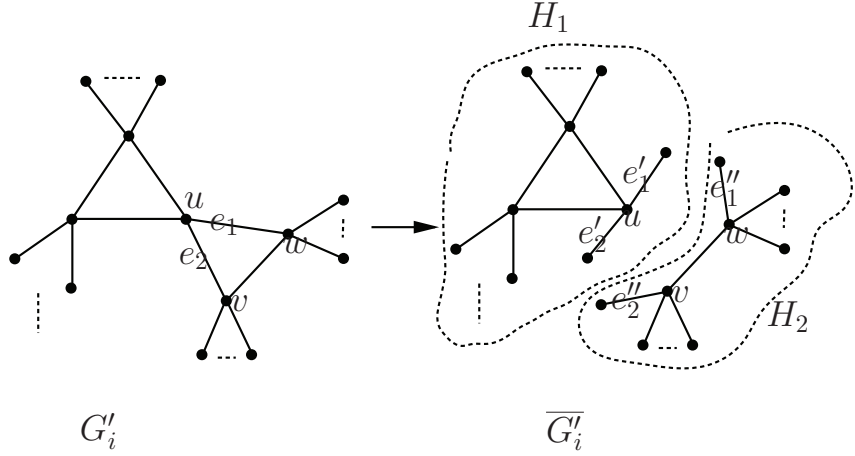}
\end{center}
\caption{$\overline{G_i'}$ is obtained from $G_i$ by doing Operation
1 to edges $e_1$ and $e_2$.}\label{figure4}
\end{figure}

For the sharpness of the upper bound, see Example \ref{exe1}.\qed
\end{pf}

We call a set of triangles $independent$ if any two of them are
vertex-disjoint. Since each single triangle is a
triangle-tree-structure, applying Theorem \ref{bsthm3}, we have the
following corollary, and for the sharpness of the upper bound, see
Example \ref{exe1}:

\begin{cor}\label{bscor1} If $G$ is a connected graph with $t'$
independent triangles, then $$rc(L(G))\leq n_2-t'.$$ Moreover, the
bound is sharp.\qed
\end{cor}

\begin{exe}\label{exe1} As shown in Figure \ref{figure5}, $G$
consists of $t$ (independent) triangles and $t-1$ edges which do not
belong to any triangles. Since $G$ has $3t$ inner vertices, by
Theorem \ref{bsthm3}(Corollary \ref{bscor1}), we know $rc(L(G))\leq
2t$; on the other hand, it is easy to show that the diameter of the
line graph $L(G)$ is $2t$, and so we have $rc(L(G))= 2t$.

\begin{figure}[!hbpt]
\begin{center}
\includegraphics[scale=1.0000000]{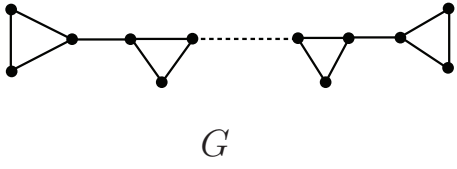}
\end{center}
\caption{Figure of Example \ref{exe1}.}\label{figure5}
\end{figure}
\end{exe}

\subsection{Another upper bound from Theorem \ref{bsthm3}}

Now we use the notation similar to that of Theorem \ref{bsthm3}. Let
$\mathcal {T}=\bigcup_{i=1}^{c}{\mathcal
{T}_i}=\bigcup_{i=1}^{c}\{T_{i,j_i}~is~ a~ triangle~ of~ $G$: 1\leq
j_i \leq t_i \}$($\sum_{i=1}^c{t_i}=t$) be a set of $t$
edge-disjoint triangles of $G$; the subgraph of $G$, $G[E(\mathcal
{T}_i)]$, induced by each $E(\mathcal {T}_i)$ is a connected
component of the subgraph $G[E(\mathcal {T})]$ which may not be a
triangle-tree-structure, that is, it may contain a
triangle-cycle-structure as a subgraph.

In $G$, for each $1\leq i\leq c$, let $G_i=G[E(\mathcal {T}_i)]$,
$V_i=V(G_i)$, $E_i=E(\mathcal {T}_i)$; $E_i^0=E(V_i)\cup
\partial(V_i)\supseteq E_i$, and $G_i^0=G[E_i^0]$. We obtain a new
graph $G'$ from $G$ by doing Operation 1 at each edge $e\in
E(V_i)\backslash E_i$ for $1\leq i\leq c$, and we denote by $G_i'$
the new subgraph (of $G'$) corresponding to $G_i^0$. Applying
Observation \ref{bsthm2} repeatedly, we have $rc(L(G))\leq
rc(L(G'))$. Now each $G_i'$ contains exactly $t_i$ triangles:
$T_{i,j_i}$ where $1\leq j_i \leq t_i$. We now obtain a new graph
$G''$ from $G'$ by doing Operation 2 to those $G_i'$s which contain
triangle-cycle-structures such that each subgraph (of $G''$) $G_i''$
corresponding to $G_i'$ contains no triangle-cycle-structure. Let
$op(G')$ be the minimum times of doing Operation 2 needed during
above procedure. Clearly, $op(G')=op(G[E(\mathcal {T})])$ (minimum
times of doing Operation 2 needed for $G[E(\mathcal {T})]$ such that
the resulting graph contains no triangle-cycle-structure). Since
there are $op(G')$ new inner vertices totally produced, and by
Observation \ref{obs2} and the discussion of Theorem \ref{bsthm3},
we have
\begin{lem}\label{bslem3.1} For any set $\mathcal {T}$ of $t$ edge-disjoint
triangles of a connected graph $G$ with $n_2$ inner vertices, we
have $$rc(L(G))\leq n_2+op(G[E(\mathcal {T})])-t.$$ \qed
\end{lem}

We know that a triangle-forest-structure of size $\ell$ contains
$2\ell+c$ (inner) vertices where $c$ is the number of components of
it. Operation 2 does not change the number of edge-disjoint
triangles, but the number of inner vertices increases 1 after we did
Operation 2 once. Then it is easy to show that after doing Operation
2 $op(G[E(\mathcal {T})])$ times, the number of inner vertices of
the new graph is $op(G[E(\mathcal {T})])$ + $n_2$=$2t+1$+$n_2'$
where $n_2'$ denotes the number of inner vertices not covered by the
original $t$ edge-disjoint triangles. So, by Lemma \ref{bslem3.1},
we have the following result and for the sharpness of the bound see
Example \ref{exe2}.

\begin{thm}\label{bsthm3.2} If $G$ is a connected graph, $\mathcal {T}$
is a set of $t$ edge-disjoint triangles that cover all but $n_2'$
inner vertices of $G$ and $c$ is the number of components of the
subgraph $G[E(\mathcal {T})]$, then
$$rc(L(G))\leq t+n_2'+c.$$ Moreover, the bound is sharp.\qed
\end{thm}

\begin{exe}\label{exe2} Let $G$ be a graph shown in Figure
\ref{figure10}. The set $\mathcal
{T}=\{u_i,v_i,u_{i+1}\}_{i=1}^{k-1}$ is a set of $k-1$ edge-disjoint
triangles, $n_2'=1$ and $c=1$. By Theorem \ref{bsthm3.2}, we have
$rc(L(G))\leq k+1$; on the other hand, it is easy to show that the
diameter of $L(G)$ is $k+1$, and so $rc(L(G))=k+1$. \qed

\begin{figure}[!hbpt]
\begin{center}
\includegraphics[scale=1.0000000]{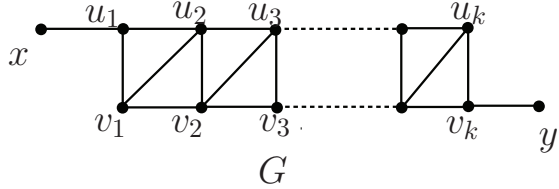}
\end{center}
\caption{Figure of Example \ref{exe2}.}\label{figure10}
\end{figure}

\end{exe}

\section{Bounds for iterated line graphs}

Recall that the $iterated~line~graph$ of a graph $G$, denoted by
$L^2(G)$, is the line graph of the graph $L(G)$. The following
corollary deduced from Theorem \ref{bslem3.1} is an upper bound of
the rainbow connection numbers of iterated line graphs of connected
cubic graphs.

\begin{cor}\label{bscor4} If $G$ is a connected cubic graph with $n$
vertices, then $$rc(L^2(G))\leq n+1.$$
\end{cor}
\begin{pf} Since $G$ is a connected cubic graph, each vertex is an
inner vertex and the clique $\langle S(v)\rangle$ in $L(G)$
corresponding to each vertex $v$ is a triangle. We know that
$\mathscr{K}=\{\langle S(v)\rangle: v\in \overline{V_2} \}=\{\langle
S(v)\rangle: v\in V \}$ is a clique decomposition of $L(G)$. Let
$\mathcal {T}=\mathscr{K}$. Then $\mathcal {T}$ is a set of $n$
edge-disjoint triangles that cover all vertices of $L(G)$ and
$L(G)=L(G)[E(\mathcal {T})]$. So $n_2'=0$ and $c=1$, and by Theorem
\ref{bsthm3.2}, the conclusion holds.\qed
\end{pf}

In graph $G$, we call a path of length $k$ a
$pendent$~$k$-$length$~$path$ if one of its end vertex has degree 1
and all inner vertices has degree 2. By definition, a pendent
$k$-length path contains a pendent $\ell$-length path($1\leq \ell
\leq k$). A pendent 1-length path is a $pendent$~$edge$.

\begin{thm} Let $G$ be a connected graph with $m$ edges and $m_1$
pendent 2-length paths. Then $rc(L^2(G))\leq m-m_1$, the equality
holds if and only if $G$ is a path of length at least 3.
\end{thm}
\begin{pf} Now $L(G)$ is a graph with $m$ vertices and $m_1$
pendent edges. Then it has $m-m_1$ inner vertices. By the discussion
before Theorem \ref{bsthm3}, we give each clique $\langle
S(v)\rangle$ in $L^2(G)$ a fresh color, where $v$ is an inner vertex
of $L(G)$. It is easy to show that this gives a rainbow
$(m-m_1)$-edge-coloring of $L^2(G)$, and so $rc(L^2(G))\leq m-m_1$.

If $G$ is a path of length at least 3, then the equality clearly
holds.

If $G$ contains a cycle, then $L(G)$ clearly contains a minimal
cycle $C: v_1,\cdots,v_\ell$. So $L^2(G)$ contains a
clique-cycle-structure of size $\ell$ which is formed by cliques
$\{\langle S(v)\rangle\}_{i=1}^\ell$. By a theorem in
\cite{Li-Sun1}, we know that this structure needs at most $\lceil
\frac{\ell+1}{2}\rceil$ colors to make sure that it is rainbow
connected. The rest part of $L^2(G)$ is formed by cliques $\{\langle
S(v)\rangle:v\in V^2\backslash \{v_i\}_{i=1}^\ell$\} where $V^2$ is
the set of inner vertices of $L(G)$. We give each $\langle
S(v)\rangle$ a fresh color, by Observation \ref{bslem2.2},
$rc(L^2(G)\leq \lceil \frac{\ell+1}{2}\rceil+(m-m_1-\ell)<m-m_1$.

If $G$ is a tree with a vertex of degree at least $3$, then $L(G)$
contains a cycle, a similar argument will show $rc(L^2(G))< m-m_1$.

So $G$ is a tree with maximum degree 2, and hence it must be a path
of length at least 3. \qed
\end{pf}

\end{document}